\documentclass{amsart}
\usepackage{latexsym,amsmath,amssymb,amsfonts,amscd,graphics}
\usepackage[mathscr]{eucal}
\usepackage[all]{xypic}

\numberwithin{equation}{section}
\theoremstyle{plain}
\newtheorem{thm}[equation]{Theorem}
\newtheorem*{mainthm}{Main Theorem}
\newtheorem{lemma}[equation]{Lemma}
\newtheorem{prop}[equation]{Proposition}
\newtheorem{corollary}[equation]{Corollary}

\newtheorem{remark}[equation]{Remark}

\newcommand{\bP}{\mathbb{P}}

\newcommand{\bQ}{\mathbb{Q}}
\newcommand{\fX}{\mathfrak{X}}
\newcommand{\sC}{\mathcal{C}}
\newcommand{\sE}{\mathcal{E}}
\newcommand{\sF}{\mathcal{F}}
\newcommand{\sK}{\mathcal{K}}
\newcommand{\sM}{\mathcal{M}}
\newcommand{\sO}{\mathcal{O}}
\newcommand{\sX}{\mathcal{X}}
\newcommand{\sD}{\mathcal{D}}
\newcommand{\Der}{\mathrm{Der}}
\newcommand{\Ext}{\mathrm{Ext}}
\newcommand{\Hilb}{\mathrm{Hilb}}
\newcommand{\Sym}{\mathrm{Sym}}
\newcommand{\PGL}{\mathrm{PGL}}
\newcommand{\Proj}{\mathrm{Proj}}
\newcommand{\gquot}{/\!\!/}
\newcommand{\spec}{\mathrm{Spec}\;}

\author{Charles Cadman and Radu Laza}
\title[Counting the hyperplane sections of a quintic]{{\bf Counting the hyperplane sections with fixed invariants of a plane quintic}\\
{\tiny Three approaches to a classical enumerative problem}} 
\address{Department of Mathematics, University of British Columbia, Vancouver, BC V6T 1Z2, Canada}
\email{cadman@math.ubc.ca}
\address{Department of Mathematics, University of Michigan, Ann Arbor, MI 48109, USA}
\email{rlaza@umich.edu}

\begin{document}
\bibliographystyle{amsplain}

\begin{abstract}
We use three different methods to count the number of lines in the plane whose intersection with a fixed general quintic has fixed cross-ratios.  We compare and contrast these methods, shedding light on some classical ideas which underly modern techniques.
\end{abstract}
\maketitle

\section*{Introduction}
\label{sec:intro}
Enumerative geometry has historically been very influential in the development of new techniques in algebraic geometry.  An example from the past decade is the persistent use of stable maps to count rational curves in Fano varieties.  In the thesis of the first author, the more recent innovation of twisted stable maps (Abramovich-Vistoli \cite{abrvis}) was used to count curves having specified tangency conditions to a hypersurface.   This paper uses the same moduli spaces to solve a different type of enumerative problem: instead of the tangency conditions,  we fix the moduli of the intersection.  This requires a type of descendent Gromov-Witten invariants where one fixes the isomorphism type of the source curve
 and  a degeneration axiom in the context of twisted stable maps.  With these tools at our disposal the enumerative problem considered here is answered quite easily. However, to prove that the answer that we obtain is enumerative and to illustrate the main ideas in a transparent way, we include also two classical approaches to the same enumerative question.

\smallskip

Concretely, we solve the following classical enumerative problem. Fix a general plane quintic curve $D$ and a configuration of $5$ points in $\bP^1$. How many lines in the plane intersect  the curve $D$ in a configuration projectively equivalent to the given one? A more precise equivalent question is to find the degree of  the generically finite rational map 
$$\phi:\check{\bP}^2\dashrightarrow \Sym^5\bP^1\gquot \PGL(2)$$
which sends a line $L$ to the moduli of the intersection $D\cap L$ (regarded as a quintuple of unordered points in $L\cong \bP^1$). Our main result is to answer to this question as follows:
\begin{mainthm}
For a general plane quintic $D$, the natural map $$\phi:\check{\bP}^2\dashrightarrow \Sym^5\bP^1\gquot \PGL(2)$$ is dominant, generically finite, and has degree $420$.
\end{mainthm}

\smallskip

The map $\phi$ appears previously in literature, and several results are known. Specifically, in general, one considers a fixed smooth degree $d\ge 3$  hypersurface $X$ in the projective space $\bP^n$. Then there exists a natural rational map $\phi: \check\bP^n\dashrightarrow \bP^{N_d}\gquot\PGL(n)$ which sends a hyperplane section $H\in \check{\bP}^n$ to the moduli of the degree $d$ hypersurface $X\cap H$ in $H\cong  \bP^{n-1}$ (with  $\bP^{N_d}\cong  \bP (H^0(\bP^{n-1},\sO_{\bP^{n-1}}(d)))$).  Beauville \cite{beauvillehyp} proved that $\phi$ is non-constant for every smooth $X$, except for some examples in positive characteristic (see also McKernan \cite{mckernan}).  In characteristic zero, one expects that the differential of $\phi$ has maximal rank at a general point, in which case one says that the hyperplane sections have maximal variation in moduli.    Harris--Mazur--Pandharipande \cite[Thm 5.1]{harrismazur} proved that for any smooth hypersurface $X$ in $\bP^n$ such that $n$ is much larger than $d$ and $k$, the variation of $k$-plane sections is maximal.  An improved bound was recently obtained by Starr \cite{starr}.  Without any bounds, but only for generic $X$, van Opstall--Veliche proved that the hyperplane sections of $X$ have maximal variation in moduli \cite[Thm 1]{opstallveliche}.  In particular, in our situation, it follows that, for a general smooth plane quintic $D$, the map $\phi$ is dominant. Thus, for dimension reasons, it makes sense to ask for its degree.

\smallskip

There are a number of problems similar to the one considered in this paper.  For example, one can ask for the degree of the map $\phi$ for the variation of hyperplane sections of cubic threefolds (together with the quintic case, these are the only two examples where the dimensions of the source and the target of $\phi$ coincide). For the cases ($d=4$ and $n=2$) and ($d=3$ and $n=3$), the fibers of $\phi$ are positive dimensional, but one can restrict to a generic pencil of hyperplanes and obtain a finite map $\bar\phi:\bP^1\to \bP^1$. It turns out that in both these cases the corresponding degree is $12$ (fact related to the theory of elliptic surfaces). For all other $n$ and $d$, the dimension of the target of $\phi$ is much larger than the dimension of the source, but one can still ask for the degree of $\phi$ onto its image.   More generally, one can produce similar problems by varying hypersurface sections instead of hyperplanes (e.g. one can count the conic sections  of a plane quartic with fixed moduli).

\smallskip

We propose three different methods for solving the enumerative problem considered here. The common idea of all three methods is to degenerate the intersection $D\cap L$ to the case when $L$ is tangent to $D$. It turns out that the difficulty of the problem is a subtle contribution to the degree coming from the flexes of $D$ (N.B.  the inflectional lines are the indeterminacy points of $\phi$).
The three methods use different techniques to handle this issue.  In the first method, the classical approach, we explicitly 
 resolve the indeterminacies of the map $\phi$ by using the stable reduction. Once this is done, we compute the degree of the map $\phi$ by comparing the self-interesection of the discriminant divisor $\Delta$ in  $\Sym^5\bP^1\gquot \PGL(2)$ to the self-intersection of its pull-back. The contribution coming from the flexes is encoded in the multiplicities of exceptional divisors occurring in the pull-back of $\Delta$.

\smallskip

The second method computes the degree by essentially lifting $\phi$ to a $\Sigma_5$-cover. Specifically, we consider the  morphism $\overline{\sM}_{0,5}^D(\bP^2,1)\to\overline{\sM}_{0,5}$, where $\overline{\sM}_{0,5}^D(\bP^2,1)$ is a closed subscheme of the space of degree 1 stable maps to $\bP^2$ where an order one contact is imposed at each marked point.  One could think of $\overline{\sM}_{0,5}^D(\bP^2,1)$ as a space of relative stable maps (see Gathmann \cite{Gathmann}).  Using deformation theory, we show that over a point of $\overline{\sM}_{0,5}$ which corresponds to a three component curve, the morphism is \'etale.  Then it is straightforward to count the number of preimages.

\smallskip

Finally, as noted in the first paragraph, the third method uses the Gromov-Witten theory of stacks to compute the degree.  While this computation very quickly produces the degree of $\phi$, we stress that the proof of enumerativity is not automatic, and depends on a result from the previous sections. 

\smallskip

The difference between the three methods is that they compute the degree of slightly different maps.  The relation between the first two is explained in section \ref{comparison}.  The point is that in the stable map approach, one orders the five intersection points between the line and quintic, so that one obtains $\Sigma_5$-covers of the source and target of $\phi$ (having the effect of killing the natural monodromy). Then, the stable map approach essentially incorporates, by construction, the stable reduction computation of the first method. 
In the third method, we use a stack of twisted stable maps $\sK$ mapping to $\overline{\sM}_{0,5}$.  
The morphism we use in the second method is essentially obtained by passing to the coarse moduli space of $\sK$.  To conclude, we can say that $\sK$ and $\overline{\sM}_{0,5}^D(\bP^2,1)$ are successive refinements of  $\check{\bP}^2$. With each refinement, the computation of  the degree of $\phi$ becomes easier, but at the price of becoming more opaque and technically demanding.  

\smallskip

We close by noting that we use the following genericity assumption in our proofs: {\it $D$ is smooth and every line in $\bP^2$ intersects $D$ in at least $3$ distinct points}.   When we say {\it a general quintic $D$}, we mean precisely this condition.  In particular, $D$ has no higher order flexes and no bitangent is also a flex.  Thus, $D$ has 120 bitangents and 45 flexes (see \cite[pp. 277--282]{griffithsharris}). On the other hand, we note that if $D$ has higher flexes, the count from the main theorem might fail. A concrete   example in this sense is the Fermat quintic, for which the degree of $\phi$ is  $150$ (see \S\ref{fermat}). A similar phenomenon (where the flexes influence the degree) was observed, in a somewhat related context, by Aluffi--Faber \cite{aluffifaber}.

\subsection*{Acknowledgements}
The first author was supported by the National Science Foundation under Grant No. 0502170.
The second author would like to thank Igor Dolgachev and Robert Friedman for some helpful discussions on the topic.

\section{First Method: The Classical Approach}
\label{sec:m1}

The first method that we propose for the computation of the degree of the map $\phi$ is based on the following observation. Let $f:X\to Y$ be a generically finite map between two smooth surfaces. Then the degree of the map $f$ can be computed by the formula
\begin{equation}\label{formuladeg}
\deg(f)=\frac{(f^*\Delta)^2}{\Delta^2},
\end{equation}
where $\Delta$ is any divisor on the surface $Y$, and the square denotes the self-intersection of the corresponding divisors. We note that the  formula is still valid when $X$ and $Y$ have at worst (finite) quotient singularities.  Surfaces with quotient singularities are $\bQ$-factorial (i.e.  every Weil divisor is $\bQ$-Cartier), and therefore the intersection and the pull-back of divisors still make sense. Essentially, the only difference to the smooth case is that the intersection numbers will be rational numbers, not necessarily integers.

We want to apply (\ref{formuladeg}) in our situation, the computation of the degree of the map $\phi:\check{\bP}^2\dashrightarrow \Sym^5\bP^1\gquot \PGL(2)$ defined for a fixed  plane quintic $D$. To do this there are two steps. First, $\phi$ is only a rational map and we need to resolve it by successively blowing-up the surface $\check{\bP}^2$ in the points of indeterminacy of $\phi$. The resulting map $\tilde \phi$ will be a morphism  of surfaces with at worst quotient singularities. The next step now is to identify a divisor $\Delta$ in  $\Sym^5\bP^1\gquot \PGL(2)$ such that we can compute its pull-back via $\tilde\phi$. Once this is done the main theorem follows directly from the formula (\ref{formuladeg}). The two steps are discussed in \S\ref{secresolution} and \S\ref{secmultiplicity} respectively. 

\subsection{Invariants of binary quintics}\label{secinvariants}
The map $\phi$ sends a line $L$ in $\bP^2$ to the moduli of the intersection $D\cap L$, regarded as a quintuple of unordered points in $\bP^1$. Here we collect a few standard facts about the moduli of quintuples as needed later. To start, we recall that a moduli space for unordered quintuples in $\bP^1$ can be constructed as the GIT quotient of $\Sym^5\bP^1\cong \bP^5$ by the natural action of  $\PGL(2)$. A quintuple is stable if no three points coincide and unstable otherwise. In particular, the indeterminacy locus of the map $\phi:\check{\bP}^2\dashrightarrow \Sym^5\bP^1\gquot \PGL(2)$ is the set of lines $L$ which are tangent to $D$ with multiplicity at least $3$. 

In the moduli space of quintuples we consider the discriminant divisor $\Delta$,   parameterizing quintuples  such that two of the five points coincide. 
The following results about the space  $\Sym^5\bP^1\gquot \PGL(2)$ and the divisor $\Delta$ are well known.

\begin{prop}
The moduli space of $5$ unordered points in $\bP^1$ is isomorphic to a weighted projective space: $\Sym^5\bP^1\gquot \PGL(2)\cong W\bP(1,2,3)$. In particular, it has only quotient singularities. Via this isomorphism the discriminant divisor $\Delta$ is a section of $\sO_{W\bP(1,2,3)}(2)$, and therefore $\Delta^2=\frac{2}{3}$. 
\end{prop}
\begin{proof}
This is a classical result of the invariant theory for binary forms (see Elliott \cite{elliott}).  Namely, the ring of invariants of binary quintics $R^G$ is generated by $4$ elements $I_4$, $I_8$, $I_{12}$ and $I_{18}$ of degrees $4$, $8$, $12$ and $18$ respectively, satisfying a relation of degree $36$. It then follows that $\Sym^5\bP^1\gquot \PGL(2)\cong \Proj (R^G)\cong W\bP(1,2,3)$ (see \cite[pg. 151-152]{dolgachevgit}).  As for the discriminant divisor, it is given by $\Delta=I_4^2-128I_8$. The proposition follows.
\end{proof}

\begin{remark}\label{stablecurves}
There exists a natural identification between the stable genus $0$ curves with $5$ marked points and stable ordered quintuples of points in $\bP^1$. Therefore, there exists a natural  morphism $\overline{\sM}_{0,5}\to  \Sym^5\bP^1\gquot \PGL(2)$ which forgets the ordering of the points. It follows that $\Sym^5\bP^1\gquot \PGL(2)$ is isomorphic to the quotient of $\overline{\sM}_{0,5}$ by the natural action of $\Sigma_5$.  In particular, we note that the pullback of $\Delta$ to $\overline{\sM}_{0,5}$ is twice the boundary divisor $\tilde{\Delta}$ in $\overline{\sM}_{0,5}$.  Since $\overline{\sM}_{0,5}$ is a degree $5$ del Pezzo with the boundary divisor $\tilde{\Delta}$ consisting of sum of ten $(-1)$-curves, one obtains $(\tilde{\Delta})^2=20$ and then $\Delta^2=4\cdot 20/120 = 2/3$.
 \end{remark}

\subsection{Resolution of the map $\phi$}\label{secresolution} Assume that $D$ is a general quintic, as defined in the introduction. The points of indeterminacy of the rational map $\phi$ are precisely the $45$ inflectional lines of $D$ (viewed as points in $\check{\bP}^2$). All these $45$ points lie on the dual curve  $\check D\subset \check \bP^2$ and, in fact, are exactly the  cusps of $\check D$. A resolution of the map $\phi$ is obtained by taking a log resolution of the dual curve at its cusps.  More precisely, we have the following result:
\begin{prop}\label{resolution}
Let $D$ be a generic plane quintic. We denote by $\tilde\bP^2$ the result of the successive blow-ups of $\check\bP^2$ in the $45$ cusps of the dual curve $\check D$  which makes the resulting exceptional divisors and the strict transform of $\check D$ meet transversally (see figure \ref{blowup}).  Then the composite map $\tilde\bP^2\to \check\bP^2\dashrightarrow \Sym^5\bP^1\gquot \PGL(2)$ extends to a morphism $\tilde\phi:\tilde\bP^2\to   \Sym^5\bP^1\gquot \PGL(2)$.
\end{prop}
\begin{proof}
The claim of the proposition is of local nature, so we can concentrate on one indeterminacy point at time. Let $o$ be such a point.  The claim of the proposition is that the following sequence of $3$ blow-ups resolves the indeterminacy at $o$. First, we blow-up $o\in \check D\subset \check \bP^2$, next we blow-up the intersection of the exceptional divisor $E_1$ with the strict transform of the dual curve $\check D$, and then we do another blow-up of the point of intersection of the new exceptional divisor $E_2$ with the strict transform of $E_1$ and of the dual curve $\check D$ (see figure \ref{blowup}).

\begin{figure}[htb]
\begin{center}
\scalebox{.40}{\includegraphics{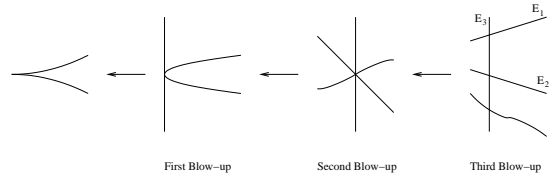}}
\end{center}
\caption{The log resolution of the dual curve}
\label{blowup}
\end{figure} 

To prove the claim we need to understand the behavior of arcs approaching $o$ under the composition with the rational map $\phi$. 

\begin{lemma}\label{lindeterm}
Let $\phi:X\dashrightarrow Y$ be a rational map between two (projective) surfaces. Assume that $X$ is smooth and $o\in X$ is an arbitrary point.  Then $\phi$ can not be extended to a morphism in a neighborhood of $o$ iff there exist two analytic arcs $f$ and $g$ from the unit disk $B$ to  $X$ with $f(0)=g(0)=o$ and such that $\phi\circ f:B^* \to Y$ and $\phi\circ g:B^* \to Y$ are well defined, but $\lim _{t\to 0} \phi(f(t))\neq \lim _{t\to 0} \phi(g(t))$. \qed
\end{lemma} 

In our situation $o\in \check \bP^2$ corresponds to an inflectional line $L_0$. An arc $f:B\to \check\bP^2$ approaching $o$ corresponds to a family of lines $(L_t)_{t\in B}$ approaching the inflectional line $L_0$. For each $t$ the quintic $D$ determines a binary quintic $D_t$ (corresponding to the intersection $D\cap L_t$). This data can be fitted in a degenerating family $(\sX,\sD)$ of $5$-pointed rational curves as follows. The family of lines can be trivialized to form the family $\sX=\bP^1\times B\to B$, where $B$ is the unit disk. Then, the binary quintics $D_t$ fit together to give a divisor $\sD$ of relative degree $5$ in $\sX\to B$. To find the limit of the arc $f$ is  equivalent to doing stable reduction (see \cite[\S3.C]{harrismorrison}) for  $(\sX,\sD)$. Namely, after a possible base change and some birational transformations affecting only the central fiber, we can transform $(\sX,\sD)\to \Delta$ into a family of stable $5$-pointed rational curves $(\sX',\sD')\to \Delta$. The new central fiber (a stable 5-pointed curve) is well defined and represents the limit $\lim _{t\to 0} \phi(f(t))$ via the identification discussed in  Remark \ref{stablecurves}.  The stable reduction procedure for a generic arc is illustrated in figure \ref{stablered} below.

\begin{figure}[htb]
\begin{center}
\scalebox{.40}{\includegraphics{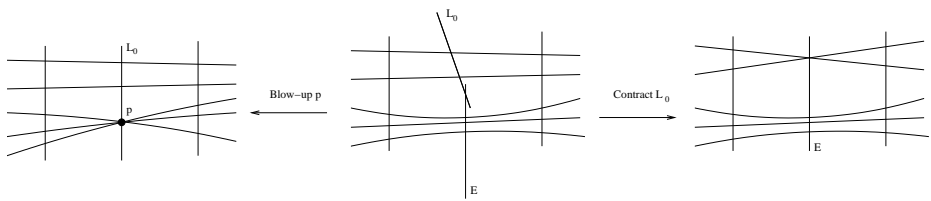}}
\end{center}
\caption{Stable reduction for a generic arc}
\label{stablered}
\end{figure} 

Note that the limit quintuple consists either of $4$ distinct points, one of which is counted with multiplicity $2$, or of $3$ distinct points, two of them counted with multiplicity $2$. In any case, the limit $\lim _{t\to 0} \phi(f(t))$ belongs to the discriminant divisor $\Delta$ in $\Sym^5\bP^1\gquot\PGL(2)$.   It is easily seen that $\Delta\cong \bP^1$, the coordinate being  the $j$-invariant of the reduced configuration ($j=\infty$ corresponds to the case of two double points). 

In coordinates, we can assume that $L_0$ is given by the equation $(x_2=0)$ and the quintic $D$ by the equation $(x_0^3(x_0-x_1)x_1+x_2f_4(x_0,x_1,x_2)=0)$ with respect to some choice of coordinates $(x_0:x_1:x_2)$ on $\bP^2$ (N.B. $(0:1:0)\in \bP^2$ is a flex of $D$). The assumption that $D$ is smooth gives that $f_4(0,1,0)\neq 0$. Without loss of generality, we can further assume $f_4(0,1,0)=1$. The family of lines $(L_t)$ have the equations $(x_2=\alpha(t)x_0+\beta(t)x_1)$ for some power series $\alpha(t)$ and $\beta(t)$ with $\alpha(0)=\beta(0)=0$. By substituting, we get a family of binary quintics:
\begin{equation}\label{famquintic}
(x_0^3(x_0-x_1)x_1+(\alpha(t)x_0+\beta(t)x_1)f_4(x_0,x_1,\alpha(t)x_0+\beta(t)x_1)=0)
\end{equation}
which specializes at $(t=0)$ to $(x_0^3(x_0-x_1)x_1=0)$. As noted above, $\sX$ is the trivial family $\sX=\bP^1\times \Delta\to \Delta$. A point on $\sX$ is denoted by $((x_0:x_1),t)$, where $t$ is an affine coordinate $t$ for the base $\Delta$ and $(x_0:x_1)$ are homogeneous coordinates for the fiber $\bP^1$. With these notations the divisor $\sD$ is given by the equation (\ref{famquintic}).

We now proceed to the explicit computation of the stable reduction:

{\bf Step 1} (Computation for the first blow-up). {\it We claim that after the blow-up of the cusp at $o$, the only indeterminacy point is the intersection point of the exceptional divisor $E_1$ and of the strict transform $\tilde D$ of the dual curve $\check D$}.  

By using lemma \ref{lindeterm}, this can be reformulated in terms of arcs as follows. Let $\alpha(t)=\alpha_0\cdot t^{n}+(\textrm{higher order terms})$ and $\beta(t)=\beta_0\cdot t^{n}+(\textrm{higher order terms})$ correspond to an arc $f:\Delta\to \check \bP^2$ (with either $\alpha_0\neq 0$ or $\beta_0\neq 0$). The claim is equivalent to saying that if $\beta_0\neq 0$, then the limit $\lim _{t\to 0} \phi(f(t))$ depends only on the fraction $\frac{\alpha_0}{\beta_0}$ and not on the arc $f$ (N.B. $(\alpha_0:\beta_0)$ is the coordinate on the exceptional divisor $E_1$). In fact,  as long as $\beta_0\neq 0$, we obtain that the limit configuration is a quintuple with $j$-invariant $0$ in the sense explained above.

After a possible base change (of type $t\to t^3$) we can assume:
\begin{eqnarray*}
\alpha(t)&=&\alpha_0\cdot t^{3k}+(\textrm{higher order terms}),\\
\beta(t)&=&\beta_0\cdot t^{3k}+(\textrm{higher order terms})
\end{eqnarray*} 
with $\beta_0\neq 0$. The equation (\ref{famquintic}) becomes:
\begin{equation}\label{famquintic2}
(x_0^3(x_0-x_1)x_1+t^{3k}(\alpha_0\cdot x_0+\beta_0\cdot x_1+t\cdot \eta(t,x_0,x_1))f_4(\dots)=0)
\end{equation}
Next, we rewrite the above equation w.r.t. the affine coordinates $x=\frac{x_0}{x_1}$ and $t$ around the point $p=((0:1),0)\in \sX$. We get:
\begin{equation}\label{famquintic3}
(x^3(x-1)+t^{3k}(\alpha_0\cdot x+\beta_0 +t\cdot \eta(t,x,1))f_4(x,1,\alpha(t)x+\beta(t))=0).
\end{equation}
After a weighted blow-up at $p$: 
$$v\cdot x=u\cdot t^k.$$ 
we obtain  a stable model, with the central fiber consisting of $2$ rational curves $L_0$ and $E$ as in figure \ref{stablered}. By contracting $L_0$ we obtain a semi-stable quintuple of points in $E\cong \bP^1$. It will consist of a double point (coming from the contraction of $L$) and three other distinct points. To compute the $j$-invariant of the reduced quadruple of points we do the following computation.  On $E\cong \bP^1$ we have the coordinates $(u:v)$. The line $L_0$ meets $E$ in the point at $\infty=(1:0)$.  The remaining $3$ points of intersection of $\sD$ with $E$ are obtained by  doing the blow-up computation, i.e.  set $x=u\cdot t^k$ (i.e. $v=1$), divide by $t^{3k}$, and then set $t=0$. We obtain
\begin{equation}
u^3+\beta_0 \cdot f_4(0,1,0)=u^3+\beta_0= 0
\end{equation}
Thus on $E$ we obtain the double point $\infty$ together with the $3$ cubic roots of the constant $-\beta_0\neq 0$. This configuration has the $j$-invariant equal to $0$.  This concludes the first step (i.e. $\beta_0\neq 0$ implies that the limit is independent of the arc and, in fact, has $j$-invariant $0$).

{\bf Step 2} (Computation for the second blow-up). We consider now the second blow-up. This corresponds to studying the limit quintuple for arcs of type $\alpha(t)=\alpha_0\cdot t^{n}+(\textrm{higher order terms})$ and $\beta(t)=\beta_0\cdot t^{m}+(\textrm{higher order terms})$ with $\alpha_0\neq 0$ and $m>n$.  We claim that if $n\le m-n$  or $\beta_0=0$ the limit does not depend on the arc, and it has invariant $j=1728$. In other words, {\it after the second blow-up the only point of indeterminacy is the intersection of the two exceptional divisors $E_1$ and $E_2$} (which is also a point on the strict transform $\tilde D$). 

The computation is similar to that from step 1. First, after a possible base change, we can assume $n=2k$. As before, we get:
\begin{equation}\label{famquintic22}
(x^3(x-1)+t^{2k}(\alpha_0\cdot x+x\cdot t\cdot \eta(x,t)+t^{m-2k}\cdot \zeta(t))f_4(x,1,\alpha(t)x+\beta(t))=0).
\end{equation}
By blowing-up $x\cdot v=t^k\cdot u$ we get
$$v^3+\alpha_0\cdot v=0$$
(under the assumption $m-2k\ge 2k>k$). As in the first step we get that the semi-stable quintuple consists of the points  $0$, $\infty$ (counted with multiplicity $2$), and $\pm \sqrt{-\alpha_0}$ on $E$. This configuration has $j$-invariant $1728$ and the claim follows.

{\bf Step 3} (Computation for the third blow-up). As in the previous step, we have to consider arcs of type $\alpha(t)=t^{n}+(\textrm{higher order terms})$ and $\beta(t)=\beta_0\cdot t^{m}+(\textrm{higher order terms})$ with $m>n$ and $\beta_0\neq 0$. Since $m\ge 2n$ was handled in the previous step, it remains to consider the case $n<m<2n$, corresponding to arcs that pass through the intersection point of $E_1\cap E_2$ after the second blow-up.  It is easily seen, by computations similar to those from steps 1 and 2, that if $3n<2m$ the limit quintuple has $j=1728$, and respectively if $3n>2m$ the limit quintuple has $j=0$ (these arcs correspond to the points of intersection of $E_3$ with $E_1$ and $E_2$ respectively). Thus, it remains only to consider arcs with $3n=2m$, i.e. $\alpha(t)=\alpha_0\cdot t^{2k}+(\textrm{higher order terms})$ and $\beta(t)=\beta_0\cdot t^{3k}+(\textrm{higher order terms})$ with $\alpha_0\neq 0$ and $\beta_0\neq 0$. A blow-up of type $x\cdot v=t^k\cdot u$ gives after simplification the equation:
$$v^3+\alpha_0\cdot v+\beta_0=0.$$
 The limit quintuple consists of the point $\infty$ (counted with multiplicity $2$) together with the roots of the previous equation. If the discriminant $4\alpha_0^3+27\beta_0^2$  vanishes, we get a configuration of $3$ distinct points, two of which are counted with multiplicity $2$. Otherwise, we obtain that the $j$-invariant for the limit is $1728 \frac{4\alpha_0^3}{4\alpha_0^3+27\beta_0^2}$, which depends only on the fraction $\frac{\alpha_0^3}{\beta_0^2}$. On the other hand, a simple computation shows that the arc given by $(\alpha(t),\beta(t))$ hits the exceptional divisor $E_3$ in the point of coordinates $(\alpha_0^3:\beta_0^2)$. Thus, {\it the map $\phi$ can be extended to $E_3$ giving a morphism $\tilde\phi: \tilde\bP^2\to \Sym^5\bP^1\gquot \PGL(2)$}. 
\end{proof}
 
In conclusion, for a general quintic (having no higher flexes), the map $\phi$ can be resolved by blowing-up three times each of the cusps of the dual curve. As a result, we obtain for each cusp $o_i$ three exceptional divisors  $E_1^{(i)}$, $E_2^{(i)}$ and $E_3^{(i)}$. The only incidences involving the exceptional divisors are that $E_3^{(i)}$ meets transversally $E_1^{(i)}$, $E_2^{(i)}$ and $\tilde D$ (the strict transform of the dual curve) as in figure \ref{blowup}. A simple computation gives the self-intersection numbers: $(E_1^{(i)})^2=-3$, $(E_2^{(i)})^2=-2$, $(E_3^{(i)})^2=-1$ and $\tilde{D}^2=130$ (N.B. the dual curve $\check D$ has degree $20$ in $\check \bP^2$).

\begin{remark}\label{corresolution}
In the proof the previous proposition we have actually shown that $E_1^{(i)}$ and $E_2^{(i)}$ are contracted via $\tilde \phi$ to points on the discriminant curve $\Delta$ and that $E_3^{(i)}\cong \bP^1$ maps isomorphically onto $\Delta\cong \bP^1$.
\end{remark}

\subsection{Computation for a general quintic}\label{secmultiplicity} For the computation of the degree of the map $\phi$ it is natural to apply the formula (\ref{formuladeg}) to the morphism $\tilde{\phi}$ and the discriminant divisor $\Delta$ in  $\Sym^5\bP^1\gquot \PGL(2)$. The essential idea is that, for geometric reasons, it is not hard to compute the pull-back $\tilde\phi^* \Delta$.  

\begin{prop}
With notations as in the previous section, we have
$$\tilde{\phi}^*\Delta=\widetilde{D}+\sum_{i=1}^{45}\left(\frac{2}{3}E_1^{(i)}+E_2^{(i)}+2E_3^{(i)}\right).$$
In particular, $\left(\tilde{\phi}^*\Delta\right)^2=280$.
\end{prop}
\begin{proof}
It is easy to see that,  set theoretically, the inverse image of $\Delta$ under $\tilde \phi$ consists of the strict transform $\tilde D$ of the dual curve and of the exceptional divisors $E_k^{(i)}$ (for $k=1,2,3$ and $i$ running over the cusps of $\tilde D$): $\tilde{\phi}^{-1}(\Delta)=\widetilde{D}\cup\bigcup_i\left( E_1^{(i)}\cup E_2^{(i)}\cup E_3^{(i)}\right)$. 
The deformation results from \S\ref{sec:deformations} (esp. Cor. \ref{discr_unramified}) say that $\phi$ is unramified at a generic point of the dual curve. Therefore, the coefficient of $\widetilde{D}$ in $\tilde\phi^*\Delta$ is precisely $1$, and we can write 
$$\tilde\phi^*\Delta=\widetilde{D}+\sum_{i}\left(a_iE_1^{(i)}+b_iE_2^{(i)}+c_iE_3^{(i)}\right)$$
 for some rational numbers $a_i$, $b_i$ and $c_i$. Recall that for a morphism $f:X\to Y$ we have the projection formula (see \cite[pg. 426, A.4]{hartshorne}):
$$f_*(f^*\beta\cdot \alpha)=\beta\cdot f_*\alpha$$
where $\alpha$ and $\beta$ are divisors on $X$ and $Y$ respectively. By applying the projection formula for  $\alpha=E_1^{(i)}$ and $\beta=\Delta$ and noting that  $f_*E_1^{(i)}=0$ ($E_1^{(i)}$ is contracted), we get 
$$f^*\Delta\cdot E_1^{(i)}=-3a_i+c_i=0 \ \Longrightarrow \ a_i=\frac{1}{3}c_i.$$
Similarly, $\alpha=E_2^{(i)}$ gives $b_i=\frac{1}{2}c_i$. Finally, by taking $\alpha=E_3^{(i)}$, we get 
$$f^*\Delta \cdot E_3^{(i)}=1+a_i+b_i-c_i=1-\frac{1}{6}c_i=\Delta\cdot\tilde\phi_* E_3^{(i)}=\Delta^2=\frac{2}{3}.$$
Therefore,  $c_i=2$ and the proposition follows.
\end{proof}

In conclusion, we obtain a proof  of the Main Theorem:
 
\begin{corollary}
The morphism $\phi$ of the introduction has degree $420$.
\end{corollary}
\begin{proof}
Clearly $\deg\phi=\deg \tilde \phi$. By (\ref{formuladeg})  and the previous proposition, we get
$$\deg \tilde\phi=\frac{(\tilde \phi^*\Delta)^2}{ \Delta^2}=\frac{3}{2}\cdot 280=420.$$
\end{proof}

\subsection{Computation for the Fermat quintic}\label{fermat}
For the computation of the degree of the map $\phi$, we made use several times of the assumption that $D$ has no higher flexes. Here, we give an example (the Fermat quintic) to show that without this genericity assumptions  the count from the main theorem  fails. 

\begin{prop}
Let $D$ be the Fermat quintic in $\bP^2$. The degree of the natural map $\phi:\check \bP^2\dashrightarrow \Sym^5\bP^1\gquot \PGL(2)$ which associates to a line $L$ the moduli of the intersection $D\cap L$  is $150$.
\end{prop}
\begin{proof}
Let $(x:y:z)$ be the coordinates on $\bP^2$. The equation of  $D$ can be taken to be  $(x^5+y^5+z^5=0)$. The map $\phi$ associates to a point $(a:b:c)\in \check \bP^2$, giving the line  $L$ of equation $ax+by+cz=0$ in $\bP^2$, the invariants of the intersection $C\cap L$. These can be computed explicitly as follows. By the linear transformation $x=\frac{x}{a}$, $y=\frac{y}{b}$ and $z=\frac{z}{c}$ the equations of the pair $(C,L)$ become
\begin{eqnarray*}
C&:& lx^5+my^5+nz^5=0,\\
L&:& x+y+z=0
\end{eqnarray*}
where $l=b^5c^5$, $m=a^5c^5$ and $n=a^5b^5$. In this form the invariants of the intersection $C\cap L$ are computed by standard formulas in the invariant theory of binary quintics. Specifically,  the invariants of the intersection $C\cap L$ are
\begin{eqnarray*}
I_4&=&(mn+nl+lm)^2-4lmn(l+m+n)=\sigma_3^2(\sigma_1^2-4\sigma_2)\\
I_8&=&(lmn)^2(mn+nl+lm)=\sigma_3^5\sigma_1\\
I_{12}&=&(lmn)^4=\sigma_3^8
\end{eqnarray*}
(cf. \cite[pg. 304]{elliott}), where $\sigma_1=a^5+b^5+c^5$, $\sigma_2=b^5c^5+a^5c^5+a^5b^5$ and $\sigma_3=a^5b^5c^5$. In conclusion, for the Fermat quintic the rational map $\phi:\check\bP^2\dashrightarrow \Sym^5\bP^1\gquot \PGL(2)\cong W\bP^2(1:2:3)$ is the composition:
\begin{eqnarray*}
 (a:b:c)\to(a^5:b^5:c^5)\to (\sigma_1:\sigma_2:\sigma_3)&\dashrightarrow& \left(\sigma_3^2(\sigma_1^2-4\sigma_2):\sigma_3^5
 \sigma_1:\sigma_3^8\right)\\
 &=&(\sigma_1^2-4\sigma_2:\sigma_1\sigma_3:\sigma_3^2).
 \end{eqnarray*}
It is easy to see that the first two maps are morphisms of degree $25$ and $6$ respectively, while the third map is a birational isomorphism of $W\bP(1:2:3)$. The claim follows.
\end{proof}

\section{Second Method: The Stable Map Approach}
\label{sec:m2}

We now demonstrate the main theorem using the moduli space $\overline{\sM}_{0,5}(\bP^2,1)$ of degree $1$ genus $0$ stable maps into $\bP^2$ with five marked points.  Note that this moduli space is actually a scheme, since none of the maps have nontrivial automorphisms.  We begin by defining a closed subscheme $$\overline{\sM}^D_{0,5}(\bP^2,1)\subseteq\overline{\sM}_{0,5}(\bP^2,1)$$ by the condition that the pullback of $D$ is the sum of the marked points.  When the source curve is reducible, we impose this condition on the degree one component.  Then we analyze the forgetful morphism $\overline{\sM}^D_{0,5}(\bP^2,1)\to\overline{\sM}_{0,5}$, whose degree is the same as that of the rational map $\phi$ in the Main Theorem.  We compute its degree by showing that it is \'etale over any point of $\overline{\sM}_{0,5}$ corresponding to a curve with three components.  It is then straightforward to count the number of preimages of such a curve.

The subscheme $\overline{\sM}^D_{0,5}(\bP^2,1)$ agrees with a space of relative stable maps $\overline{\sM}_{\alpha}^D(\bP^2,1)$ as defined by Andreas Gathmann \cite{Gathmann}, with $\alpha=(1,1,1,1,1)$.  He defined them for genus $0$ maps relative to a very ample divisor.  In positive genus, an algebro-geometric definition was given by Jun Li \cite{JunLi}.  In the case we are interested in, it is not hard to construct the space by hand and work out the deformation theory we need directly.  The fact that any degree $1$ genus $0$ stable map has a ``main component'' which does not map into the quintic $D$ means that we can compare $D$ with the marked points on this main component and see whether they are equal.  This condition cuts out a locus of $\overline{\sM}_{0,5}(\bP^2,1)$ scheme-theoretically.

In fact, this point of view gives us a functor from schemes to sets which sends a scheme $S$ the set of degree $1$, genus $0$, $5$-marked stable maps to $\bP^1$ over $S$ where the condition ``$D$ = marked points'' is verified after collapsing all degree $0$ components in the fibers of the curve over $S$.  From this, we can compute the tangent space to our subscheme by applying this functor to the spectrum of $\mathbb{C}[x]/(x^2)$.  This is implicit in our calculation of the deformation spaces and helps us show that the map $\overline{\sM}^D_{0,5}(\bP^2,1)\to\overline{\sM}_{0,5}$ is unramified over the ``most degenerate'' curves in $\overline{\sM}_{0,5}$.

It is worth mentioning that $\overline{\sM}^D_{0,5}(\bP^2,1)$ is the closure inside of $\overline{\sM}_{0,5}(\bP^2,1)$ of the set of lines in $\bP^2$ transverse to $D$ with the intersection points marked.

\subsection{Definition of $\overline{\sM}^D_{0,5}(\bP^2,1)$}
\label{sec:def}

Let $\sC_{0,n}(\bP^2,1)$ be the universal curve over $\overline{\sM}_{0,n}(\bP^2,1)$.  Recall that
\begin{equation}
\label{can_isom}
\sC_{0,n}(\bP^2,1)\cong\overline{\sM}_{0,n+1}(\bP^2,1)
\end{equation}
in such a way that the projection to $\overline{\sM}_{0,n}(\bP^2,1)$ is the map forgetting the last marked point and the universal morphism $f_n:\sC_{0,n}(\bP^2,1)\to\bP^2$ is the evaluation map $e_{n+1}$.  We have a commutative diagram.
$$\xymatrix{
\sC_{0,5}(\bP^2,1) \ar[r] \ar[d] & \sC_{0,0}(\bP^2,1) \ar[d] \\
\overline{\sM}_{0,5}(\bP^2,1) \ar@/^/[u]^{s_i} \ar[r] & \overline{\sM}_{0,0}(\bP^2,1)}$$
The horizontal arrows forget the first $5$ marked points and $s_i$ is the $i$th section, $1\le i\le 5$.  Note that $\overline{\sM}_{0,0}(\bP^2,1)$ is just $\check{\bP}^2$, and $\sC_{0,0}(\bP^2,1)$ is the flag variety of points and lines in $\bP^2$.  Let $$\tilde{\sC}_{0,5}(\bP^2,1)=\overline{\sM}_{0,5}(\bP^2,1)\times_{\overline{\sM}_{0,0}(\bP^2,1)}\sC_{0,0}(\bP^2,1),$$ with first and second projections $\pi_1,\pi_2$.  The natural morphism $c:\sC_{0,5}(\bP^2,1)\to\tilde{\sC}_{0,5}(\bP^2,1)$ contracts all the degree $0$ components of the fibers.  Let $\tilde{s_i}=c\circ s_i$.  Since the fibers of $\pi_1$ are smooth, it follows that the image of $s_i$ is an effective Cartier divisor (see Cadman \cite[5.1]{cadman1}).  Let $$\Sigma\subseteq\tilde{\sC}_{0,5}(\bP^2,1)$$ be the sum of these five divisors.

By the Remark following \cite[5.2]{cadman1}, the closed subscheme $\pi_2^{-1}f_0^{-1}D$ of $\tilde{\sC}_{0,5}(\bP^2,1)$ is an effective Cartier divisor.  We define $\overline{\sM}^D_{0,5}(\bP^2,1)$ to be the locus where this divisor agrees with $\Sigma$.  To make this precise, let $\Hilb^5_{\pi_1}$ be the relative Hilbert scheme parametrizing length five subschemes of the fibers of $\pi_1$.  Then $\Sigma$ and $\pi_1^*f_0^*D$ define sections of the projection $\Hilb^5_{\pi_1}\to\overline{\sM}_{0,5}(\bP^2,1)$.  We define $\overline{\sM}^D_{0,5}(\bP^2,1)$ to be the fiber product.
$$\xymatrix{
\overline{\sM}^D_{0,5}(\bP^2,1) \ar[r] \ar[d] \ar@{}[rd]|{\Box} & \overline{\sM}_{0,5}(\bP^2,1) \ar[d]^{\pi_2^*f_0^*D} \\
\overline{\sM}_{0,5}(\bP^2,1) \ar[r]^{\Sigma} & \Hilb^5_{\pi_1}}$$

Using this definition, we claim that $\overline{\sM}^D_{0,5}(\bP^2,1)$ is a closed subscheme of $\overline{\sM}_{0,5}(\bP^2,1)$ whose dimension at each point is at least two.  First note that $\tilde{\sC}_{0,5}(\bP^2,1)$ is the projectivization of a rank $2$ vector bundle over $\overline{\sM}_{0,5}(\bP^2,1)$.  Therefore $\Hilb^5_{\pi_1}$ is the projectivization of the fifth symmetric power of this bundle, hence smooth of relative dimension 5 over $\overline{\sM}_{0,5}(\bP^2,1)$.  It follows that any section of $\Hilb^5_{\pi_1}$ is a regular embedding (EGAIV \cite[17.12.1]{grothendieck}), from which the claim follows.

With a little more work, we can show that $\overline{\sM}^D_{0,5}(\bP^2,1)$ is purely two-dimensional.  Under the genericity assumption on $D$, every line intersects $D$ in at least three points.  It follows that the map $$\overline{\sM}^D_{0,5}(\bP^2,1)\to\overline{\sM}_{0,0}(\bP^2,1)=\check{\bP}^2$$ is finite away from the flexes of $D$.  Over the flexes, the fibers are at most one dimensional.  Therefore, the dimension of $\overline{\sM}^D_{0,5}(\bP^2,1)$ cannot be larger than two.

\subsection{Deformations fixing a divisor}
\label{sec_deform_divisor}

Here we prove an easy lemma, which is essential for the deformation-obstruction theory of relative stable maps (see Graber-Vakil \cite[2.8]{GraberVakil}).  This lemma is well-known and often left as an exercise, but we include a proof for completeness.

\begin{lemma}
\label{lemma_deformations}
Let $f:X\to Y$ be a morphism of smooth varieties over a field $k$, and let $D\subseteq Y$ be a smooth divisor such that $f^{-1}D$ is also a divisor.  Let $\sE$ be the sheaf of logarithmic vector fields on $Y$ relative to $D$ (often denoted $T_Y(-\mathrm{log}\, D)$).  Then the first order deformations of $f$ which fix $X$, $Y$, and $f^{-1}D$ correspond naturally to $H^0(X,f^*\sE)$.
\end{lemma}

\begin{proof}
Recall that $\sE$ is the kernel of the natural morphism $T_Y\to\sO_D(D)$, and is locally free of rank $2$.  Since $f^{-1}D$ is a divisor, it follows that $L_1f^*\sO_D(D)=0$, so the following sequence is exact.
$$0\to f^*\sE\to f^*T_Y\to f^*\sO_D(D)\to 0$$
Suppose we have a first order deformation of $f:X\to Y$ fixing $X$ and $Y$.  This corresponds to an element $\xi\in H^0(X,f^*T_Y)$.  We need to show that $\xi$ preserves $f^*D$ if and only if the image of $\xi$ in $H^0(X,f^*\sO_D(D))$ is zero.  This is local, so assume $Y=\spec R$ and $X=\spec S$.  Let $\psi=f^{\#}:R\to S$, and let $\alpha\in R$ be an equation for $D$.  Then $\psi(\alpha)\neq 0$ by assumption.

The morphism $H^0(X,f^*T_Y)\to H^0(X,f^*\sO_D(D))$ is locally represented by a map 
$$\xymatrix@R=0pt{
\Der_k(R,S) \ar[r] & (\psi(\alpha)^{-1}\cdot S)/S \\
\rho \ar@{|->}[r] & \frac{\rho(\alpha)}{\psi(\alpha)}.}$$
The first order deformation $R\to S[\epsilon]/(\epsilon^2)$ corresponding to $\xi\in\Der_k(R,S)$ is given by $r\mapsto\psi(r)+\xi(r)\epsilon$.

This preserves $f^*D$ if and only if there is a unit $u\in S[\epsilon]/(\epsilon^2)$ such that $u\psi(\alpha)=\psi(\alpha)+\xi(\alpha)\epsilon$.  Since $S$ is a domain and every element of the form $1+v\epsilon$ is a unit, this happens if and only if $\xi(\alpha)$ is in the ideal generated by $\psi(\alpha)$.  This is equivalent to $\xi$ mapping to zero in $H^0(X,f^*\sO_D(D))$.
\end{proof}

\subsection{Infinitesimal analysis}\label{sec:deformations}

Let $\sC=(f:C\to\bP^2,x_1,\ldots,x_5)$ represent a closed point of $\overline{\sM}^D_{0,5}(\bP^2,1)$, where $x_i$ are the marked points.  Let $C_0\subset C$ be the degree $1$ component.  Since $C_0$ intersects $D$ at least three times, and each intersection point is either a node or marked point of $C_0$, it follows that $C_0$ has at least three special points.  It is well-known that $H^1(C,f^*T\bP^2)=0$ for any genus $0$ stable map.  Thus we have an exact sequence.
$$0\to H^0(C,f^*T\bP^2) \to T_{\sC}\overline{\sM}_{0,5}(\bP^2,1) \to \Ext^1(\Omega_C(\sum x_i),\sO_C) \to 0$$
Note that $\Ext^1(\Omega_C(\sum x_i),\sO_C)$ is the tangent space to $\overline{\sM}_{0,5}$ at $(C,x_1,\ldots,x_5)$.

We investigate whether the morphism $$T_{\sC}\overline{\sM}^D_{0,5}(\bP^2,1)\to \Ext^1(\Omega_C(\sum x_i),\sO_C)$$ has a nonzero kernel.  An element of the kernel corresponds to a first order deformation of $\sC$ which fixes both $C$ and $f_0^*D$, where $f_0=f\vert_{C_0}$.  By Lemma \ref{lemma_deformations}, the kernel is identified with $H^0(C_0,f_0^*\sE)$, where $\sE$ is the kernel of $T_{\bP^2}\to\sO_D(D)$.  Since $f_0^*\sE$ is a rank $2$ degree $-2$ vector bundle on $C_0\cong\bP^1$, we have the following.

\begin{prop}
\label{prop_etale}
The morphism $\overline{\sM}^D_{0,5}(\bP^2,1)\to\overline{\sM}_{0,5}$ is \'etale at $\sC$ if and only if there is an injection of sheaves $\sO_{C_0}(-1)\to f_0^*\sE$ whose cokernel is locally free.
\end{prop}

\begin{proof}
The latter condition is equivalent to $f_0^*\sE$ being isomorphic to $\sO_{C_0}(-1)^{\oplus 2}$, which is equivalent to the vanishing of $H^0(C_0,f_0^*\sE)$.
\end{proof}

\begin{prop}
\label{prop_jump_lines}
If $C_0$ is either a flex or a bitangent of $D$, then there is an injection of sheaves $\sO_{C_0}(-1)\to f_0^*\sE$ whose cokernel is locally free.
\end{prop}

\begin{proof}
First we do a local computation to compute the fiber of $f_0^*\sE$ at a point of $C_0$.  Suppose $D$ is defined at $(0,0)$ by the (analytic) equation $y=g(x)$, where $g^{(k)}(0)\ne 0$ and $g^{(m)}(0)=0$ for $m<k$.  Let $C_0$ be the line $y=0$.  Then $C_0$ has a $k$th order contact with $D$.

If $v = P(x)\partial/\partial x + Q(x)\partial/\partial y$ is a local vector field along $C_0$, then $v$ is a local section of $f_0^*\sE$ if and only if $v(y-g(x))$ belongs to the ideal generated by $g(x)$ (by definition of $\sE$).  Up to a normalization factor, this means that $x^{k-1}P(x)$ is congruent to $Q(x)$ modulo $x^k$.  It now follows that the fiber of $f_0^*\sE$ at $x=0$ has as a basis the images of $\partial/\partial x + x^{k-1}\partial/\partial y$ and $x\partial/\partial x$.

Let $\sF\subset T_{C_0}$ be the subsheaf consisting of sections which vanish at every intersection point with $D$.  By hypothesis, there are exactly three intersection points, so $\sF\cong\sO_{C_0}(-1)$.  It follows from the local computation that we have an injection $\sF\to f_0^*\sE$ whose cokernel is locally free.
\end{proof}

\begin{corollary}\label{discr_unramified}
The morphism $\overline{\sM}_{0,5}^D(\bP^2,1)\to\overline{\sM}_{0,5}$ is unramified over a general point of the boundary divisor.  The same holds for the morphism $\phi$ in the Main Theorem over the discriminant divisor.
\end{corollary}

\begin{figure}[htb]
\begin{center}
\scalebox{.60}{\includegraphics{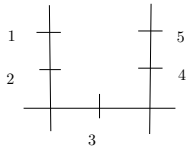}}
\end{center}
\caption{A curve in $\overline{M}_{0,5}$}
\label{fig_curve}
\end{figure}

\begin{proof}
Let $C\in \overline{\sM}_{0,5}$ be the curve illustrated in figure \ref{fig_curve}.  Referring back to the construction of $\overline{\sM}_{0,5}^D(\bP^2,1)$, it is clear that any stable map in $\overline{\sM}_{0,5}^D(\bP^2,1)$ whose underlying curve is $C$ corresponds either to a bitangent (with the middle component having degree $1$) or to a flex (with one of the other two components having degree $1$).  By Propositions \ref{prop_etale} and \ref{prop_jump_lines}, the morphism is unramified at any of these points.  Since any component of the boundary divisor contains a curve in the $\Sigma_5$-orbit of $C$, it follows that the branch divisor does not contain any such component.

For the second part of the corollary, note that at a general point of the boundary divisor in $\overline{\sM}_{0,5}$, the isotropy subgroup of $\Sigma_5$ fixing both that point and any of its preimages is the $\mu_2$ which interchanges the two marked points which are isolated on a component.  This means that the property of being \'etale is preserved after taking quotients by $\Sigma_5$.
\end{proof}

\begin{remark}
\emph{From the results of this section, we can say a little more about the ramification divisor of $\phi$, namely that it is a curve of degree $12$ in $\check{\bP}^2$.  For this we use a result of Barth \cite{barth} which characterizes the set of jumping lines of a vector bundle in projective space.  From it, one can deduce that the set of jumping lines of $\sE$ (lines $L$ for which the restriction $\sE\vert_L$ does not split into $\sO_L(-1)^{\oplus 2}$) form a curve of degree $12$.  One must be careful with the isotropy groups to make sure that $\phi$ does not have extra (divisorial) ramification which is not present in the 
$\Sigma_5$-cover, but this is not hard.}
\end{remark}

\subsection{Degree of the morphism}\label{degreemet2}

\begin{thm}
\label{thm_degree}
The morphism $\overline{\sM}_{0,5}^D(\bP^2,1)\to\overline{\sM}_{0,5}$ has degree $420$.
\end{thm}

\begin{proof}
Let $C$ be the curve in figure \ref{fig_curve}.  Any preimage of $C$ intersects $D$ in exactly three points, so is either a bitangent or a flex.  By Propositions \ref{prop_etale} and \ref{prop_jump_lines}, the morphism is \'etale over $C$, so the degree of the morphism is the number of preimages of $C$.  Each bitangent counts twice since the components mapping to the tangent points can be interchanged.  Each flex counts four times, since either of two components can map with positive degree and the marked points on that component can be interchanged.  So the degree is $2\cdot 120+4\cdot 45=420$.
\end{proof}

We now compare this with the Main Theorem.  There is a forgetful morphism $$\overline{\sM}_{0,5}(\bP^2,1)\to\overline{\sM}_{0,0}(\bP^2,1)\cong \check{\bP}^2$$ which is fixed by the $\Sigma_5$ action permuting the marked points, so it induces a morphism $\overline{\sM}_{0,5}^D/\Sigma_5\to\check{\bP}^2$.  So we have the following commutative diagram: 
\begin{equation}\label{diagram1}
\xymatrix{
\overline{\sM}_{0,5}^D\ar[r]\ar[d] & \overline{\sM}_{0,5}^D/\Sigma_5 \ar[d] \ar[r] \ar[dr]^{\bar{\phi}} & \check{\bP}^2 \ar@{-->}[d]^-{\phi} \\
\overline{\sM}_{0,5} \ar[r] & \overline{\sM}_{0,5}/\Sigma_5 \ar[r]^-{\sim} & \Sym^5\bP^1\gquot\bP GL(2)}
\end{equation}
where $\bar{\phi}$ is the composition of the induced morphism $\overline{\sM}_{0,5}^D/\Sigma_5 \to \overline{\sM}_{0,5}/\Sigma_5$ with the isomorphism $\overline{\sM}_{0,5}/\Sigma_5\cong\Sym^5\bP^1\gquot\bP GL(2)$.  We also use $\overline{\sM}_{0,5}^D$ to refer to $\overline{\sM}_{0,5}^D(\bP^1,1)$. From this it follows immediately the Main Theorem: 

\begin{corollary}\label{cor_main_thm_2}
The morphism $\phi$ of the introduction has degree $420$.
\end{corollary}

\begin{proof}
In the above diagram the action of $\Sigma_5$ is generically free on both $\overline{\sM}_{0,5}^D$ and $\overline{\sM}_{0,5}$.  Moreover, the morphism $\overline{\sM}_{0,5}^D/\Sigma_5\to\check{\bP}^2$ is birational, so the degree of $\phi$ is the same as that of the morphism on the left.
\end{proof}

\subsection{Comparison of the two methods}\label{comparison}
By construction, a stable pointed curve (in the sense of $\overline{\sM}_{g,n}$) is the result of applying stable reduction to a family degenerating curves. Thus, roughly speaking, the space $\overline{\sM}_{0,5}^D$ is a systematic packing of the stable reduction computations used in the proof of \ref{resolution}. It follows that there should be a close relation between $\overline{\sM}_{0,5}^D$ and the space $\tilde \bP^2$. Indeed, this is the case: {\it $\tilde \bP^2$ is a minimal resolution of the singularities $\overline{\sM}_{0,5}^D/\Sigma_5$}, giving the following  extension of the diagram (\ref{diagram1}) (N.B. factoring by the natural $\Sigma_5$-action simply amounts to passing from ordered to unordered quintuples).

\begin{equation}\label{diagram2}
\xymatrix{
&& \tilde{\bP}^2 \ar[dd]_(.3){\tilde{\phi}} \ar[ld]_-{\epsilon} \ar[rd] & \\
\overline{\sM}_{0,5}^D\ar[r]\ar[d]  &\overline{\sM}_{0,5}^D/\Sigma_5 \ar[dr]^-{\bar{\phi}} \ar[rr] \ar[d]& & \check{\bP}^2 \ar@{-->}[dl]^-{\phi} \\
 \overline{\sM}_{0,5} \ar[r] & \overline{\sM}_{0,5}/\Sigma_5 \ar[r]^-{\sim} & \Sym^5\bP^1\gquot\bP GL(2) & }
 \end{equation}

\begin{figure}[htb]
\begin{center}
\scalebox{.60}{\includegraphics{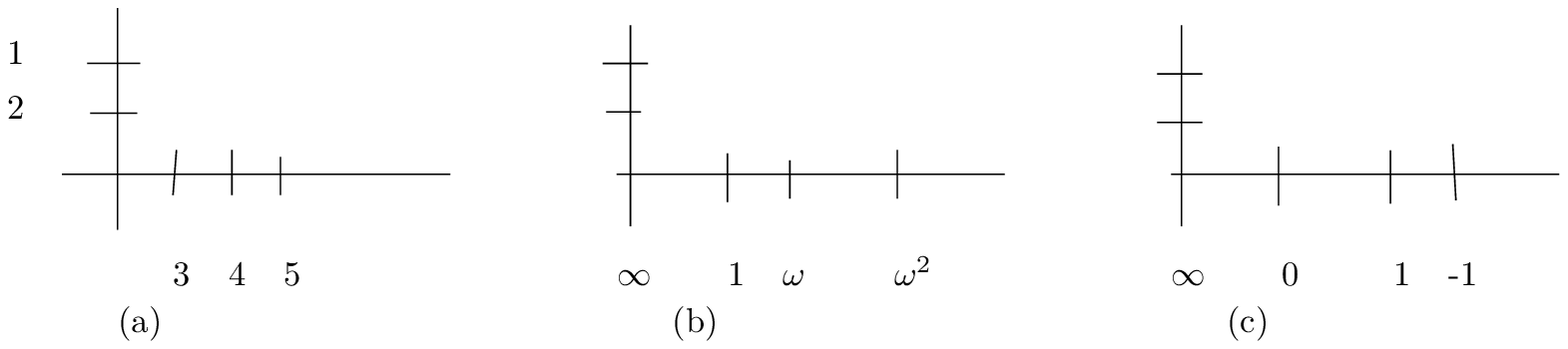}}
\end{center}
\caption{}
\label{fig_curves}
\end{figure}

To see this one should note that  each flex of $D$ gives  $2$ quotient singularities, of types $\frac{1}{3}(1,1)$ and $\frac{1}{2}(1,1)$ respectively, for  $\overline{\sM}_{0,5}^D/\Sigma_5$. Namely, consider the points in $\overline{\sM}_{0,5}^D$ mapping to a point in $\check \bP^2$ that corresponds to  an inflectional line to $D$. These points are  a union of rational curves in $\overline{\sM}_{0,5}^D$ parameterizing the  maps from a  two (or  three) component curve (see figure \ref{fig_curves}a) onto an inflectional line such that the 345 component is sent to the flex point. Forgetting the labels, we obtain for each flex, exactly one curve $E^{(i)}$ in $\overline{\sM}_{0,5}^D/\Sigma_5$ that is contracted to a point in $\check \bP^2$. At a generic point of $E^{(i)}$ the space $\overline{\sM}_{0,5}^D/\Sigma_5$ is smooth, but at special points of  $E^{(i)}$ (with large stabilizers)  $\overline{\sM}_{0,5}^D/\Sigma_5$ becomes singular. Specifically, the  two singular points mentioned above correspond to the case  illustrated in figures \ref{fig_curves}b and \ref{fig_curves}c: there are a $\mu_3$ and a $\mu_2$ stabilizer respectively that give rise to  quotient singularities for $\overline{\sM}_{0,5}^D/\Sigma_5$.  For example, in figure \ref{fig_curves}b, the  $\mu_3$ which cyclically permutes markings 3,4, and 5 does not change the curve, because it is the same as multiplication by $\omega$ on $\bP^1$ (with $\omega=e^{2\pi i/3}$).

To conclude, one can see that (see section \ref{secresolution}) that there is a natural identification  of $E_3^{(i)}$  to $E^{(i)}$, and that (at least set theoretically) $E_1^{(i)}$ and $E_2^{(i)}$ correspond to the two singular points of $\overline{\sM}_{0,5}^D/\Sigma_5$ that lie on $E^{(i)}$ (see also remark \ref{corresolution}). With a little more work, one can show that indeed there exists a morphism $\epsilon: \tilde\bP^2 \to \overline{\sM}_{0,5}^D/\Sigma_5$ mapping $E_3^{(i)}$ isomorphically onto $E^{(i)}$ and contracting $E_1^{(i)}$ and $E_2^{(i)}$. Since $\tilde\bP^2$ is smooth, we conclude that $\tilde\bP^2$ is nothing other than a resolution of $\overline{\sM}_{0,5}^D/\Sigma_5$ (which is even minimal).


\section{Third Method: Gromov--Witten for Stacks}
\label{sec:m3}

Finally, we sketch a method of computing the degree of $\phi$ which relies on the degeneration axiom for \emph{twisted} Gromov-Witten invariants.  A statement and proof of this axiom can be found in Abramovich, Graber, and Vistoli \cite[Prop. 5.3.1]{abrgravis}.

Recall that in their most general form, Gromov-Witten invariants are multilinear maps whose inputs are cohomology classes in the target space, and whose output is a cohomology class in the moduli space of marked curves (see Manin \cite[III.5]{manin}).  Here we use the word `cohomology' to mean the Chow group modulo numerical equivalence.  The degeneration axiom says what happens when the output is intersected with the boundary of $\overline{\sM}_{g,n}$ (see \cite[III.5.2.iii]{manin}).

Our target is the stack $\fX:=\bP^2_{D,r}$ for some integer $r\ge 4$.  This stack is locally the quotient of a cyclic $r$-sheeted cover of $\bP^2$ ramified along $D$ by the Galois action.  Such a cover only exists locally, but they glue to a global stack (see Cadman \cite{cadman1} for details).  The inputs for the Gromov-Witten invariants of a stack are cohomology classes on its inertia stack.  For simplicity, we'll use the cohomology of the inertia stack's coarse moduli space instead.  This has $r$ components, one of which is $\bP^2$ and the other $r-1$ of which are isomorphic to $D$.  We use the notation $\fX_0$ for $\bP^2$ and $\fX_i$ for the remaining components, $1\le i\le r-1$.  Let $\alpha_i$ be the fundamental class of $\fX_i$, and let $\beta_i$ be the class of a point on $\fX_i$.

These $r$ components should be thought of in the following way.  Let $\sK_{0,5}(\fX,1)$ be the stack of $5$-marked twisted stable maps to $\fX$ of degree $1$ and genus $0$.  Given such a map $F:\sC\to\fX$, the $j$-th evaluation map sends this map to the pair $(F(x_j),\rho_{F,j})$, where $x_j\in\sC$ is the $j$th marked point, and $\rho_{F,j}$ is the reduction mod $r$ of the order of contact between the coarse curve of $\sC$ and the divisor $D$ at the marked point $x_j$.  If $\rho_{F,j}$ is positive, then $F(x_j)$ is a point of $D$, so the pair $(F(x_j),\rho_{F,j})$ is a point of $\fX_{\rho_{F,j}}$.  Otherwise, it is a point of $\fX_0$.  This leads to the idea that using the class $\alpha_i$ as an input in a Gromov-Witten invariant is the same as imposing an $i$th order contact with the curve $D$ at an arbitrary point.  To use the class $\beta_i$ is the same as imposing it at a fixed general point of $D$.  For more details see Cadman--Chen \cite{cadmanchen}.

This leads one to speculate that the degree of the morphism $\phi$ from the introduction should equal the Gromov-Witten invariant $I_1^{\fX}(\alpha_1^5)$.  The notation means that we look at degree $1$ (genus $0$) stable maps into $\fX$, pull back the class $\alpha_1$ once for each of five evaluation maps, and push forward the resulting class to $\overline{\sM}_{0,5}$.  This yields a degree $0$ cohomology class, which we identify with an element of $\bQ$.

\begin{thm}
The Gromov-Witten invariant $I_1^{\fX}(\alpha_1^5)$ equals $420$.
\end{thm}

\begin{proof}
The degeneration axiom allows us to compute the following.  Note that it is important that $r\ge 4$.

\begin{align*}
I_1(\alpha_1^5) &= rI_0(\alpha_1^2\beta_{r-2})I_1(\alpha_1^3\alpha_2)+rI_1(\alpha_1^2\alpha_3)I_0(\alpha_1^3\beta_{r-3}) \\
I_1(\alpha_1^3\alpha_2) &= rI_0(\alpha_1^2\beta_{r-2})I_1(\alpha_1\alpha_2^2)+rI_1(\alpha_1^2\alpha_3)I_0(\alpha_1\alpha_2\beta_{r-3}) \\
I_0(\alpha_1^3\beta_{r-3}) &= rI_0(\alpha_1^2\beta_{r-2})I_0(\alpha_1\alpha_2\beta_{r-3})
\end{align*}

The following can be computed directly.
\begin{align*}
I_0(\alpha_1^2\beta_{r-2}) &= 1/r \\
I_0(\alpha_1\alpha_2\beta_{r-3}) &= 1/r \\
I_1(\alpha_1^2\alpha_3) &= 2\cdot 45 \\
I_1(\alpha_1\alpha_2^2) &= 2\cdot 120
\end{align*}

Combining, we obtain
$$I_1(\alpha_1^5)=I_1(\alpha_1\alpha_2^2)+2I_1(\alpha_1^2\alpha_3) = 420.$$
\end{proof}

As a consequence of the previous theorem we obtain again the Main Theorem:
\begin{corollary}
The morphism $\phi$ of the introduction has degree $420$.
\end{corollary}

\begin{proof}
Let $\sK\subseteq\sK_{0,5}(\fX,1)$ be the intersection of the preimages of $\fX_1$ under all five evaluation maps.  The theorem shows that the morphism $\sK\to\overline{\sM}_{0,5}$ has degree $420$.  It suffices to see that the morphism $\sK\to\overline{\sM}_{0,5}(\bP^2,1)$ maps $\sK$ birationally onto $\overline{\sM}_{0,5}^D(\bP^2,1)$, since then one recovers Theorem \ref{thm_degree}.  First, from the fact that $r\ge 4$ and that $D\subset\bP^2$ is general, one can deduce that for every map $F:\sC\to\fX$ in $\sK$, every preimage of a point of $D$ is a twisted point (untwisted points only occur when $r$ divides the order of contact).  Moreover, it is not hard to see then that $\sK\to\overline{\sM}_{0,5}^D(\bP^2,1)$ is a bijection.

We claim that it is an isomorphism away from the boundary.  Suppose $F:\sC\to\fX$ is a twisted stable map with $\sC$ irreducible.  Let $\sE$ be the pullback to $\fX$ of the sheaf of logarithmic vector fields on $\bP^2$ relative to $D$ (which is the same as those on $\fX$ relative to the gerbe over $D$).  By Cadman \cite[3.4.3]{cadman2}, we have isomorphisms $$H^i(\sC,F^*\sE)\to H^i(\sC,F^*T_{\fX}).$$  Since the deformation theory of $\sC$ coincides with that of its (marked) source curve, the claim follows from Lemma \ref{lemma_deformations}.  This shows that $\sK$ maps birationally onto $\overline{\sM}_{0,5}^D(\bP^2,1)$, which completes the proof.
\end{proof}

\bibliography{degreferences}

\providecommand{\bysame}{\leavevmode\hbox to3em{\hrulefill}\thinspace}
\providecommand{\MR}{\relax\ifhmode\unskip\space\fi MR }
\providecommand{\MRhref}[2]{%
  \href{http://www.ams.org/mathscinet-getitem?mr=#1}{#2}
}
\providecommand{\href}[2]{#2}
\begin{thebibliography}{10}

\bibitem{abrgravis}
D.~Abramovich, T.~Graber, and A.~Vistoli, \emph{Gromov--{W}itten theory of
  {D}eligne--{M}umford stacks}, arXiv:math.AG/0603151 (2006), 47 pp.

\bibitem{abrvis}
D.~Abramovich and A.~Vistoli, \emph{Compactifying the space of stable maps}, J.
  Amer. Math. Soc. \textbf{15} (2002), no.~1, 27--75 (electronic).

\bibitem{aluffifaber}
P.~Aluffi and C.~Faber, \emph{Linear orbits of smooth plane curves}, J.
  Algebraic Geom. \textbf{2} (1993), no.~1, 155--184.

\bibitem{barth}
W.~Barth, \emph{Some properties of stable rank-{$2$} vector bundles on {${\bf
  P}\sb{n}$}}, Math. Ann. \textbf{226} (1977), no.~2, 125--150.

\bibitem{beauvillehyp}
A.~Beauville, \emph{Sur les hypersurfaces dont les sections hyperplanes sont
  \`a module constant}, The Grothendieck Festschrift, Vol.\ I, Progr. Math.,
  vol.~86, Birkh\"auser Boston, Boston, MA, 1990, pp.~121--133.

\bibitem{cadman2}
C.~Cadman, \emph{Gromov-{W}itten invariants of {$\mathbb{P}\sp 2$}-stacks},
  Compos. Math. \textbf{143} (2007), no.~2, 495--514.

\bibitem{cadman1}
\bysame, \emph{Using stacks to impose tangency conditions on curves}, Amer. J.
  Math. \textbf{129} (2007), no.~2, 405--427.

\bibitem{cadmanchen}
C.~Cadman and L.~Chen, \emph{Enumeration of rational plane curves tangent to a
  smooth cubic}, arXiv:math.AG/0701406 (2007), 27 pp.

\bibitem{dolgachevgit}
I.~V. Dolgachev, \emph{Lectures on invariant theory}, London Mathematical
  Society Lecture Note Series, vol. 296, Cambridge University Press, Cambridge,
  2003.

\bibitem{elliott}
E.~B. Elliott, \emph{An introduction to the algebra of quantics}, second ed.,
  Oxford at the Clarendon Press, 1913.

\bibitem{Gathmann}
A.~Gathmann, \emph{Absolute and relative {G}romov-{W}itten invariants of very
  ample hypersurfaces}, Duke Math. J. \textbf{115} (2002), no.~2, 171--203.

\bibitem{GraberVakil}
T.~Graber and R.~Vakil, \emph{Relative virtual localization and vanishing of
  tautological classes on moduli spaces of curves}, Duke Math. J. \textbf{130}
  (2005), no.~1, 1--37.

\bibitem{griffithsharris}
P.~Griffiths and J.~Harris, \emph{Principles of algebraic geometry}, Wiley
  Classics Library, John Wiley \& Sons Inc., New York, 1994.

\bibitem{grothendieck}
A.~Grothendieck, \emph{\'{E}l\'ements de g\'eom\'etrie alg\'ebrique. {IV}.
  \'{E}tude locale des sch\'emas et des morphismes de sch\'emas {IV}}, Inst.
  Hautes \'Etudes Sci. Publ. Math. (1967), no.~32, 361.

\bibitem{harrismazur}
J.~Harris, B.~Mazur, and R.~Pandharipande, \emph{Hypersurfaces of low degree},
  Duke Math. J. \textbf{95} (1998), no.~1, 125--160.

\bibitem{harrismorrison}
J.~Harris and I.~Morrison, \emph{Moduli of curves}, Graduate Texts in
  Mathematics, vol. 187, Springer-Verlag, New York, 1998.

\bibitem{hartshorne}
R.~Hartshorne, \emph{Algebraic geometry}, Springer-Verlag, New York, 1977,
  Graduate Texts in Mathematics, No. 52.

\bibitem{JunLi}
J.~Li, \emph{Stable morphisms to singular schemes and relative stable
  morphisms}, J. Differential Geom. \textbf{57} (2001), no.~3, 509--578.

\bibitem{manin}
Y.~I. Manin, \emph{Frobenius manifolds, quantum cohomology, and moduli spaces},
  American Mathematical Society Colloquium Publications, vol.~47, American
  Mathematical Society, Providence, RI, 1999.

\bibitem{mckernan}
J.~McKernan, \emph{Varieties with isomorphic or birational hyperplane
  sections}, Internat. J. Math. \textbf{4} (1993), no.~1, 113--125.

\bibitem{starr}
J.~Starr, \emph{Fano varieties and linear sections of hypersurfaces},
  arXiv:math.AG/0607133 (2006), 8 pp.

\bibitem{opstallveliche}
M.~A. van Opstall and R.~Veliche, \emph{Variation of hyperplane sections},
  arXiv:math.AG/0602137 (2006), 4 pp.

\end{thebibliography}

\end{document}